\documentclass[sn-mathphys]{sn-jnl}

\newif\ifnarrowmode
\narrowmodetrue

\usepackage{bm}

\graphicspath{{fig/}}
\DeclareGraphicsExtensions{.mps,.pdf,.eps,.png}

\DeclareMathOperator{\rank}{rank}

\DeclareMathOperator{\trace}{tr}

\newcommand{\fro}{\mathsf F}

\newcommand*{\trans}{^{\top}}
\newcommand*{\herm}{^*}

\newcommand*{\iherm}{^{-*}}

\newcommand*{\me}{\mathrm e}
\newcommand*{\mi}{\mathrm i}

\def\adots{\mathinner{\mkern2mu\raise1pt\hbox{.}\mkern2mu
    \raise4pt\hbox{.}\mkern2mu\raise7pt\hbox{.}\mkern1mu}}

\newcommand*{\macheps}{\bm u}

\newcommand*{\tol}{\mathtt{tol}}

\newcommand*{\single}{\mathtt{lower}}
\newcommand*{\double}{\mathtt{working}}

\DeclareMathOperator{\fl}{f{}l}

\theoremstyle{plain}
\newtheorem{theorem}{Theorem}
\newtheorem{lemma}[theorem]{Lemma}

\theoremstyle{remark}
\newtheorem*{remark}{Remark}

\let\REQUIRE=\Require
\let\ENSURE=\Ensure
\let\STATE=\State
\let\FOR=\For
\let\ENDFOR=\EndFor
\let\IF=\If
\let\ENDIF=\EndIf

\renewcommand{\qedhere}{}

\newcommand\MSC[1]{\pacs[MSC Classification]{#1}}

\jyear{2023}%

\begin{document}

\title[Mixed Precision LOBPCG]{A mixed precision LOBPCG algorithm}


\author[1]{\fnm{Daniel} \sur{Kressner}}\email{daniel.kressner@epfl.ch}
\equalcont{These authors contributed equally to this work.}

\author[2]{\fnm{Yuxin} \sur{Ma}}\email{yxma18@fudan.edu.cn}
\equalcont{These authors contributed equally to this work.}

\author[3,4]{\fnm{Meiyue} \sur{Shao}}\email{myshao@fudan.edu.cn}
\equalcont{These authors contributed equally to this work.}

\affil[1]{\orgdiv{Institute of Mathematics},
\orgname{EPFL},
\orgaddress{\city{Lausanne}, \postcode{CH-1015}, \country{Switzerland}}}

\affil[2]{\orgdiv{School of Mathematical Sciences},
\orgname{Fudan University},
\orgaddress{\city{Shanghai}, \postcode{200433}, \country{China}}}

\affil[3]{\orgdiv{School of Data Science},
\orgname{Fudan University},
\orgaddress{\city{Shanghai}, \postcode{200433}, \country{China}}}

\affil[4]{\orgdiv{MOE Key Laboratory for Computational Physical Sciences},
\orgname{Fudan University},
\orgaddress{\city{Shanghai}, \postcode{200433}, \country{China}}}


\abstract{%
The locally optimal block preconditioned conjugate gradient (LOBPCG) algorithm
is a popular approach for computing a few smallest eigenvalues and the
corresponding eigenvectors of a large Hermitian positive definite matrix \(A\).
In this work, we propose a mixed precision variant of LOBPCG that uses a
(sparse) Cholesky factorization of \(A\) computed in reduced precision as the
preconditioner.
To further enhance performance, a mixed precision orthogonalization strategy
is proposed.
To analyze the impact of reducing precision in the preconditioner on
performance, we carry out a rounding error and convergence analysis of PINVIT,
a simplified variant of LOBPCG.
Our theoretical results predict and our numerical experiments confirm that the
impact on convergence remains marginal.
In practice, our mixed precision LOBPCG algorithm typically reduces the
computation time by a factor of \(1.4\)--\(2.0\) on both CPUs and GPUs.
}

\keywords{%
Symmetric eigenvalue problem, LOBPCG algorithm, mixed precision algorithm
}

\MSC{%
65F15, 65F50
}

\maketitle

\section{Introduction}
\label{sec:introduction}

Given a large Hermitian positive definite matrix
\(A\in\mathbb{C}^{n\times n}\), this work considers the computation of the
\(k\) smallest eigenvalues $0 < \lambda_1 \le \cdots \le \lambda_k$ and the
corresponding eigenvectors $x_1$, $\ldots$, $x_k$ satisfying
\[
AX=X\Lambda,
\]
where \(X = [x_1,\ldots,x_k]\) and \(\Lambda\) is diagonal with diagonal
entries $\lambda_1$, $\ldots$, $\lambda_k$.
This problem is often encountered in many applications, such as PDE and
optimization problem, electronic structure calculations and machine learning;
see, for example, \cite{BGHRCV2010, K2017, S2011}.

When a good preconditioner $T$ for \(A\) is available, the preconditioned
inverse iteration (PINVIT) from~\cite{N2001-1} is a good candidate for solving
such eigenvalue problems.
For \(k=1\), PINVIT takes the form
\[
x_{i+1} = x_i - T\bigl(Ax_i-\rho(x_i)x_i\bigr),
\]
for some starting vector \(x_0\).
Here, \(\rho(x) = (x\herm Ax)/(x\herm x)\) denotes the Rayleigh quotient,
which is also used to approximate the eigenvalue at each iteration.
Note that PINVIT with the ``ideal'' preconditioner \(T=A^{-1}\) becomes
equivalent to inverse iteration.
When computing several ($k>1$) smallest eigenpairs, one chooses a starting
matrix \(X_0\in \mathbb{C}^{n\times m}\) (\(m\geq k\)) with orthonormal columns
and one step of the block version of PINVIT~\cite{N2002} takes the form
\[
\Tilde X_{i+1} = X_i - T(AX_i-X_i\Theta_i),
\]
where \(\Theta_i = X_i\herm AX_i\).
The next iterate $X_{i+1}$ is obtained from orthonormalizing the columns
of~$\Tilde X_{i+1}$ by, e.g., a QR factorization.
Under mild conditions, linear convergence of PINVIT is proven
in~\cite{AKNOZ2017},  with a convergence rate depending on the quality of the
preconditioner \(T\).
The locally optimal block preconditioned conjugate gradient (LOBPCG)
method~\cite{K2001} aims at accelerating the convergence of PINVIT by choosing
the next iterate optimally from a $3m$-dimensional subspace that contains the
current as well as the previous iterate and the preconditioned residual;
see Section~\ref{sec:preliminary} for more details.
LOBPCG converges at least as fast as PINVIT and often significantly faster.

Executing an algorithm in reduced (single) precision on, e.g., a GPU, can be
significantly faster than executing it in default working (double) precision.
On the other hand, critical applications may require eigenvalues and
eigenvectors computed to an accuracy warranted by working precision.
In such a scenario the use of mixed precision algorithms can be beneficial;
see~\cite{AABCCD2021,HM2022} for an overview.
For example, Carson and Higham~\cite{CH2018} proposed a general framework for
large-scale mixed precision linear system solvers based on iterative
refinement.
It is highlighted that a mixed precision algorithm can be twice as fast as a
traditional linear system solver by computing the most expensive part---LU
factorization---in reduced precision.
For eigenvalue problems, mixed precision algorithms have recently been
proposed for computing \emph{all} eigenvalues and eigenvectors of a dense
matrix.
This includes the Newton-like iterative refinement methods for
symmetric~\cite{OA2018,OA2019,OA2020} and nonsymmetric~\cite{BKS2022}
eigenvalue problems, as well as a mixed precision one-sided Jacobi SVD
algorithm~\cite{GMS2022}.
If only a few eigenvalues and eigenvectors are of interest, one could combine
mixed precision with classical iterative refinement~\cite{D1982} for
eigenvalue problems, which solves linear systems with the shifted matrix
$A-\hat \lambda_i I$ in order to correct an approximation $\hat \lambda_i$ of
the $i$th eigenvalue.
The need for solving several differently shifted linear systems makes such an
approach rather expensive.

In this work, we propose mixed precision PINVIT and LOBPCG algorithms that use
a (sparse) Cholesky factorization of $A$ computed in reduced precision as
preconditioner.
This reduces the cost of accurately computing eigenvalues and eigenvectors in
significantly compared to inverse iteration, which requires to carry out the
Cholesky factorization in working precision.
On the theoretical side, we carry out a rounding error analysis of PINVIT,
which predicts that reducing precision in the preconditioner usually only has
a marginal impact on convergence.
On the experimental side, we demonstrate for sparse matrices that our mixed
precision LOBPCG algorithm results in up to \(1.43\times\) speedup on a CPU
and \(1.67\times\) speedup on a GPU.
For dense matrices, the speedups are \(1.67\times\) on a CPU and
\(2.00\times\) on a GPU.

The rest of this paper is organized as follows.
In Section~\ref{sec:preliminary}, we explain the basic ideas of LOBPCG
algorithm.
Then in Section~\ref{sec:algorithm}, we propose our mixed precision algorithms
and the details of the implementation.
The analysis is shown in Section~\ref{sec:convergence} and numerical
experiments are presented in Section~\ref{sec:experiments} to show the
efficiency of our mixed precision LOBPCG algorithm.

\section{LOBPCG algorithm}
\label{sec:preliminary}

In this section, we explain the basic idea of the LOBPCG algorithm
from~\cite{K2001}.
For $k = 1$, LOBPCG can be derived from the preconditioned conjugate gradient
(PCG) method.
PCG applied to the (singular) linear system \((A-\lambda_1I)x = 0\) with
preconditioner $T$ and  initial guess~\(x_0\) is a three-term recurrence of
the form
\ifnarrowmode
\begin{align*}
x_{i+1} &= x_i+\alpha_i T(A-\lambda_1I)x_i + \beta_i(x_i-x_{i-1})\\
&= (1+\beta_i)x_i + (-\beta_i)x_{i-1} + \alpha_i T(A-\lambda_1I)x_i,
\end{align*}
\else
\[
x_{i+1} = x_i+\alpha_i T(A-\lambda_1I)x_i + \beta_i(x_i-x_{i-1})
= (1+\beta_i)x_i + (-\beta_i)x_{i-1} + \alpha_i T(A-\lambda_1I)x_i,
\]
\fi
where \(\alpha_i\), \(\beta_i\) are chosen to minimize
\(
x_{i+1}\herm (A-\lambda_1I)x_{i+1}.
\)
As the smallest eigenvalue $\lambda_1$ is usually unknown, it needs to be
replaced by an approximation, the Rayleigh quotient \(\rho(x_i)\), leading to
the basic form of LOBPCG:
\[
x_{i+1} = \alpha^{(i)}_1x_i+\alpha^{(i)}_2x_{i-1}
+\alpha^{(i)}_3T(Ax_i-\rho(x_i)x_i),
\]
where \(\alpha^{(i)}_1\), \(\alpha^{(i)}_2\), and \(\alpha^{(i)}_3\) are
chosen to minimize \(\rho(x_{i+1})\).
Note that, unlike PCG, LOBPCG is not a Krylov subspace method in the usual
sense because \(\rho(x_i)\) is different in each iteration.

For \(k > 1\), LOBPCG takes an initial guess \(X_0\in \mathbb{C}^{n\times m}\)
with $m\geq k$, and produces iterates of the form
\[
X_{i+1} = X_iC_1^{(i)} + X_{i-1}C_2^{(i)} + W_iC_3^{(i)}
=\begin{bmatrix} X_i & X_{i-1} & W_i \end{bmatrix}
\begin{bmatrix} C_1^{(i)} \\ C_2^{(i)} \\ C_3^{(i)} \end{bmatrix}
=: S_iC_i,
\]
where \(W_i = T(AX_i-X_i\Theta_i)\) with \(\Theta_i = X_i\herm AX_i\).
The $3m\times m$ matrix $C_i$ is chosen to minimize
\begin{equation}
\label{eq:mintr}
\min_{X_{i+1}\herm X_{i+1} = I}\trace(X_{i+1}\herm AX_{i+1})
= \min_{C_i\herm S_i\herm S_iC_i = I}\trace(C_i\herm S_i\herm AS_iC_i),
\end{equation}
where $\trace(\cdot)$ denotes the trace of a matrix.
By the Rayleigh--Ritz method, a solution $C_i$ of~\eqref{eq:mintr} is obtained
from the eigenvectors belonging to the $m$ smallest eigenvalues of the
generalized eigenvalue problem \(S_i^* A S_i y = \lambda S_i^* S_i y\);
see~\cite[Section~8.7.2]{GV2013} for numerical algorithms.

Let us stress that the actual implementation of LOBPCG is quite
different~\cite{DSYG2018} due to the numerical instability caused by the
ill-conditioning of $S_i$.
In practice \([X_i,X_{i-1}]\) can be orthogonalized by an improved
Hetmaniuk--Lehoucq trick~\cite[Section~4.2]{DSYG2018}, and then the remaining
block, \(W_i\), also needs to be orthogonalized carefully.

\section{Mixed precision algorithms}
\label{sec:algorithm}
In this section, we derive a mixed precision LOBPCG algorithm.
For this purpose, we consider two precisions: a working precision and a
lower/reduced precision, e.g., IEEE double and single precisions.
The input and output data of our algorithms are always stored in working
precision.
The functions \(\single(\cdot)\) and \(\double(\cdot)\) are used to convert
working precision data into lower precision and vice versa.

\subsection{Lower precision preconditioning}
\label{subsec:preconditioning}

The application of the preconditioner $T$ usually consumes a considerable
fraction of the computational expense of PINVIT and LOBPCG.
This suggests to implement the application of $T$ in lower precision.
In most cases, we expect that this only has a small impact on convergence.
While a more detailed analysis will be provided in
Section~\ref{sec:convergence}, the existing convergence analysis of PINVIT
already provides a good intuition.

By~\cite[Theorem~2.1]{AKNOZ2017}, PINVIT with $k=1$ converges to the smallest
eigenvalue and eigenvector when \(\gamma:=\lVert I-A^{1/2}TA^{1/2}\rVert_2<1\)
and additional mild conditions are satisfied.
Asymptotically, the convergence is linear with a rate that is bounded by
$\gamma + (1-\gamma) \lambda_1/\lambda_2$.
If $T$ is perturbed by rounding error in lower precision one effectively
applies a preconditioner $T_E$, which remains close to $T$.
In turn, the convergence is now determined by
$\lVert I-A^{1/2}T_EA^{1/2}\rVert_2$, which remains close to~$\gamma$.
Unless $\gamma$ is very close to $1$ we thus expect that replacing $T$ by
$T_E$ does not affect convergence significantly.
These considerations lead to Algorithm~\ref{alg:mppinv}, PINVIT with a lower
precision preconditioner.

\begin{algorithm}[!tb]
\caption{Mixed precision PINVIT algorithm}
\label{alg:mppinv}
\begin{algorithmic}[1]
\REQUIRE
A Hermitian positive definite matrix \(A\in\mathbb C^{n\times n}\);
an initial approximate eigenvectors \(X_0\in\mathbb C^{n\times m}\);
the number of desired eigenpairs \(k\leq m\);
the maximum number of iterations \(\mathtt{maxit}\);
a function \(f_T(\cdot)\) to apply preconditioner \(T\) (in lower precision).
\ENSURE
The diagonal matrix \(\Theta\in\mathbb R^{k\times k}\) contains the computed
smallest eigenvalues, and \(X\in\mathbb C^{n\times k}\) contains the
corresponding computed eigenvectors satisfying \(AX = X\Theta\).

\STATE \(\Tilde X\gets X_0\).
\FOR{\(i = 1\), \(2\), \(\dotsc\), \(\mathtt{maxit}\)}
    \STATE \(X\gets Q\) where \(Q\) satisfies \(\Tilde X=QR\).
    \STATE \(\Theta\gets X\herm AX\).
    \STATE Compute the residual \(R = AX-X\Theta\).
    \IF{\(k\) smallest eigenpairs have converged}
        \STATE{Return \(X\gets X(:, 1:k)\) and
        \(\Theta \gets \Theta(1:k, 1:k)\).}
    \ENDIF
    \STATE Compute \(W_\single\gets f_T(\single(R))\) in a lower precision.
    \STATE \(W\gets \double(W_\single)\).
    \STATE \(\Tilde X\gets X - W\).
\ENDFOR
\end{algorithmic}
\end{algorithm}

\subsection{A mixed precision orthogonalization procedure}
In both PINVIT and LOBPCG, we need to produce an orthogonal basis of the
searching subspace in each iteration.
Moreover, orthogonalization plays an important role to ensure numerical
stability for the LOBPCG algorithm~\cite{DSYG2018,HL2006}.
We need to perform the orthogonalization procedure as accurately as possible.
However, orthogonalization is often quite expensive in practice.
Therefore it is desirable to make use of a lower precision to accelerate this
procedure.

There are mainly two existing mixed precision algorithms for computing the QR
factorization.
The algorithm proposed in~\cite{YTD2015} uses higher precision to compute the
inner product to enhance the numerical stability of Cholesky-QR algorithm.
The drawback is that this algorithm can be much slower than the standard
Cholesky-QR algorithm if higher precision arithmetic lacks hardware support.
To improve the performance, a mixed precision block Gram--Schmidt
orthogonalization algorithm was proposed in~\cite{YTKDB2015}.
For both algorithms the orthogonality of the output depends linearly on the
condition number of the input.

We propose another mixed precision approach for orthogonalization.
We first use Householder-QR to factorize
\(\single(W_i) = Q_\single R_\single\) in lower precision.
Then \(\double(R_\single)\) is used as a preconditioner---we apply Cholesky-QR
to the preconditioned matrix \(W_i\cdot\double(R_\single^{-1})\) to refine the
orthogonality.
Under mild assumptions \(W_i\cdot\double(R_\single^{-1})\) is reasonably
well-conditioned, so that the Cholesky-QR algorithm is sufficiently accurate.
This mixed precision QR factorization algorithm is summarized in
Algorithm~\ref{alg:mqr}.

\begin{algorithm}[!tb]
\caption{Mixed precision QR factorization algorithm}
\label{alg:mqr}
\begin{algorithmic}[1]
\REQUIRE
A matrix \(A\in\mathbb C^{n\times m}\) with \(\rank(A)=m\).
\ENSURE
A matrix \(Q\in \mathbb C^{n\times m}\) and an upper triangular matrix
\(R\in\mathbb C^{m\times m}\) satisfying \(A=QR\) and \(Q\herm Q=I_m\).

\STATE Compute the QR factorization of \(A\) in lower precision, i.e.,
\(\single(A) = Q_\single R_\single\).
\STATE Compute \(V\gets A\cdot\double(R_\single^{-1})\) by solving an upper
triangular linear system.
\STATE Compute Cholesky factorization of \(V\herm V\) such that
\(V\herm V = LL\herm\).
\STATE Compute \(Q\gets VL\iherm\) by solving an upper triangular linear
system.
\end{algorithmic}
\end{algorithm}

\subsection{A mixed precision LOBPCG algorithm}
In addition to preconditioning and orthogonalization, the application of $A$
and other parts of PINVIT and LOBPCG may also constitute nonnegligible
expenses, depending on the specific setting.
Carrying out these parts in lower precision bears the risk of limiting the
attainable accuracy to lower precision.
However, very often it is still possible to further exploit lower precision
arithmetic.

As PINVIT and LOBPCG converge linearly in general, we can break the
computation in two stages as follows.
In the first stage we can first perform all computations in lower precision
to produce an approximate solution in lower precision.
Then in the second stage we switch back to the working precision while using
the approximate solution as an initial guess and applying lower precision
preconditioning.
In this manner we are able to obtain a satisfactory solution in working
precision by making use of lower precision arithmetic as much as possible.

In summary, we compute a good initial guess in lower precision, and then
refine the solution using the LOBPCG algorithm in working precision.
Lower precision are exploited in both preconditioning and orthogonalization in
the LOBPCG algorithm.
The resulting mixed precision LOBPCG algorithm is summarized in
Algorithm~\ref{alg:mlobpcg}.

\begin{algorithm}[!tb]
\caption{Mixed precision LOBPCG algorithm}
\label{alg:mlobpcg}
\begin{algorithmic}[1]
\REQUIRE
A Hermitian positive definite matrix \(A\in\mathbb C^{n\times n}\);
an initial approximate eigenvectors \(X_0\in\mathbb C^{n\times m}\);
the number of desired eigenpairs \(k\leq m\);
the maximum number of iterations \(\mathtt{maxit}\);
a function \(f_T(\cdot)\) to apply preconditioner \(T\) (in lower precision).
\ENSURE
The diagonal matrix \(\Theta\in\mathbb C^{k\times k}\) contains the computed smallest eigenvalues and \(X\in\mathbb C^{n\times k}\) contains the corresponding computed eigenvectors satisfying \(AX = X\Theta\).

\STATE Compute \(X\) by a lower precision LOBPCG algorithm with the initial guess \(\single(X_0)\).
\STATE \(P\gets[~]\), \(\Theta\gets X\herm AX\).

\FOR{\(i = 1, 2, \dotsc, \mathtt{maxit}\)}
  \STATE Compute residual \(R\gets AX-X\Theta\).
  \STATE Determine the number of convergence eigenpairs \(n_c\).
  \IF{\(n_c\geq k\)}
    \STATE Return \(X \gets X(:, 1:k)\) and \(\Theta\gets \Theta(1:k, 1:k)\).
  \ENDIF
  \label{line:w-start}
  \STATE Compute \(W_\single\gets f_T(\single(R))\) in a lower precision.
  \STATE \(W\gets \double(W_\single)\).
  \label{line:w-end}
  \STATE Orthogonalize \(W\) against \([X,P]\) twice.
  \STATE Factorize \(W=QR\) and then replace \(W\) by \(Q\).
  \STATE \(A_p\gets S\herm AS\) with \(S = [X, P, W]\).
  \STATE Solve eigenvalue problem of \(A_p\) to obtain
  \(A_pC = CD\).
  \label{line:hl}
  \STATE Orthogonalize \(C(1:m, m+1:n)\herm\) to obtain an orthogonal matrix
  \(V\).
  \STATE \(X\gets SC(:, 1:m)\), \(P\gets SC(:, m+1:2m)V\) and \(\Theta\gets D(1:m, 1:m)\).
  \label{line:hl-end}
\ENDFOR
\end{algorithmic}
\end{algorithm}

\section{Convergence in finite-precision arithmetic}
\label{sec:convergence}

In our experiments, we observe that rounding error does not significantly
affect the convergence of Algorithms~\ref{alg:mppinv} and~\ref{alg:mlobpcg}
until an accuracy on the level of \emph{working} precision is reached.
To gain theoretical insights on this observation, we study the effect of
rounding error on PINVIT for $k=1$:
\begin{equation}
\label{eq:iterate-pinvit-2}
x_{i+1} = x_i - T\bigl(Ax_i - \rho(x_i)x_i\bigr).
\end{equation}
For simplicity, we consider real matrices, that is,
\(A \in \mathbb R^{n\times n}\) is positive definite with eigenvalues
\(0<\lambda_1<\lambda_2\leq \dotsb\leq \lambda_n\).
Moreover, we assume that \(n^{-1}\) is far larger than the unit roundoff, even
in reduced precision.

In analyzing the effect of rounding error on~\eqref{eq:iterate-pinvit-2}, we
assume that the computed matrix--vector product $\fl(Ax_i)$ satisfies the
backward error
\begin{equation}
\label{eq:er-Ax-ge}
\fl(Ax_i) = (A + \Delta A)x_i
\qquad\text{with}\qquad
\lVert\Delta A\rVert_2\leq \epsilon_{A}\lVert A\rVert_2,
\end{equation}
for some symmetric $\Delta A$ (depending on $x_i$).
When carrying out standard matrix--vector multiplication with a dense or
sparse matrix $A$ then Lemma~6.6 in \cite{H2002} states
that~\eqref{eq:er-Ax-ge} holds with
\[
\epsilon_A=\sqrt{n}\gamma_n^h
\qquad\text{with}\qquad
\gamma_n^h = \frac{n\macheps_h}{1-n\macheps_h},
\]
where $\macheps_h$ denotes the unit roundoff in working precision.

\begin{lemma}
\label{lem:r}
Let \(\hat{r}_i\) denote the result of evaluating \(r_i:=Ax_i - \rho(x_i)x_i\)
in working precision.
Assuming that~(\ref{eq:er-Ax-ge}) holds, there exist a symmetric matrix
\(F\in\mathbb{R}^{n\times n}\) and a diagonal matrix
\(E\in\mathbb{R}^{n\times n}\) such that
\[
\hat{r}_i = (I+E)(r_i+Fx_i),
\]
where \(\lVert E\rVert_2\leq \macheps_h\) and
\(\lVert F\rVert_2\leq \epsilon_r\lVert A\rVert_2\) with
\[
\epsilon_r
= \bigl(\gamma_n^h+\epsilon_A+\gamma_n^h\epsilon_A+(n+1)\macheps_h\bigr)
\frac{1+\macheps_h}{1-2n\macheps_h}+\epsilon_A+\macheps_h.
\]
\end{lemma}
\begin{proof}
We first analyze the rounding error when forming \(\rho(x_i)\).
From~\cite[Equation~(3.5)]{H2002} and~\eqref{eq:er-Ax-ge}, we obtain
\[
\bigl\lvert \fl(x_i\trans Ax_i) - x_i\trans \fl(Ax_i)\bigr\rvert
\leq \gamma_n^h\lVert x_i\rVert_2\lVert \fl(Ax_i)\rVert_2
\leq \gamma_n^h(1+\epsilon_A)\lVert A\rVert_2\lVert x_i\rVert_2^2.
\]
Thus, we have
\begin{align*}
\bigl\lvert\fl(x_i\trans Ax_i) - x_i\trans Ax_i\bigr\rvert
\leq{}&\bigl\lvert\fl(x_i\trans Ax_i) - x_i\trans \fl(Ax_i)\bigr\rvert
+\bigl\lvert x_i\trans\fl(Ax_i)-x_i\trans Ax_i\bigr\rvert\\
\leq{}&\gamma_n^h(1+\epsilon_A)\lVert A\rVert_2\lVert x_i\rVert_2^2
+ \lVert x_i\trans \Delta A x_i\rVert_2\\
\leq{}&\bigl(\gamma_n^h(1+\epsilon_A)+\epsilon_A\bigr)
\lVert A\rVert_2\lVert x_i\rVert_2^2.
\end{align*}
Combined with \(\fl(x_i\trans x_i) = x_i\trans x_i(1+\delta_1)\) for
$\lvert\delta_1\rvert \le \gamma_n^h$, this implies for
\(\rho(x_i) = x_i\trans Ax_i/(x_i\trans x_i)\) that there is
\(\lvert\delta_2\rvert\leq \macheps_h\) such that
\begin{align}
\bigl\lvert \fl(\rho(x_i))-\rho(x_i)\bigr\rvert
&= \left\lvert \frac{\fl(x_i\trans Ax_i)}
{x_i\trans x_i(1+\delta_1)}(1+\delta_2)
-\frac{x_i\trans Ax_i}{x_i\trans x_i(1+\delta_1)}
(1+\delta_2+\delta_1-\delta_2)\right\rvert \nonumber\\
&\leq \left\lvert
\frac{\bigl(\fl(x_i\trans Ax_i)-x_i\trans Ax_i\bigr)(1+\delta_2)}
{x_i\trans x_i(1+\delta_1)}\right\rvert
+\left\lvert\frac{x_i\trans Ax_i(\delta_2-\delta_1)}
{x_i\trans x_i(1+\delta_1)}\right\rvert \nonumber\\
&\leq \frac{\bigl\lvert \fl(x_i\trans Ax_i)-x_i\trans Ax_i\bigr\rvert}
{x_i\trans x_i} \left\lvert\frac{1+\delta_2}{1+\delta_1}\right\rvert
+\rho(x_i)\left\lvert\frac{\delta_2-\delta_1}{1+\delta_1}\right\rvert
\nonumber \\
&\leq \left(\gamma_n^h(1+\epsilon_A)+\epsilon_A\right)\frac{1}{1-2n\macheps_h}
\lVert A\rVert_2+\frac{(1+n)\macheps_h}{1-2n\macheps_h}\rho(x_i)\nonumber\\
&\leq \bigl(\gamma_n^h(1+\epsilon_A)+\epsilon_A+(n+1)\macheps_h\bigr)
\frac{\lVert A\rVert_2}{1-2n\macheps_h}.
\label{eq:er-rho}
\end{align}

The vector subtraction and scaling when forming \(r_i=Ax_i - \rho(x_i)x_i\)
yield two diagonal matrices \(E\) and \(E_1\) such that
\begin{align*}
\hat{r}_i&=(I+E)\bigl((A+\Delta A)x_i-(I+E_1)\fl(\rho(x_i)) x_i \bigr)\\
&=(I+E)( r_i + F x_i),\\
F&:= \Delta A-\fl(\rho(x_i)) E_1 - \bigl(\fl(\rho(x_i))-\rho(x_i)\bigr)I.
\end{align*}
where \(\lVert E\rVert_2\leq \macheps_h\) and
\(\lVert E_1\rVert_2\leq \macheps_h\).
Combined with~\eqref{eq:er-rho}, this concludes the proof because
\begin{align*}
\lVert F\rVert_2
&\leq \macheps_h\lVert A\rVert_2+
(1+\macheps_h)\bigl(\gamma_n^h(1+\epsilon_A)+\epsilon_A+(n+1)\macheps_h\bigr)
\frac{1}{1-2n\macheps_h}\lVert A\rVert_2+\epsilon_A\lVert A\rVert_2\\
&\leq \Bigl( \bigl(\gamma_n^h+\epsilon_A+\gamma_n^h\epsilon_A+
+(n+1)\macheps_h\bigr)
\frac{1+\macheps_h}{1-2n\macheps_h}+\epsilon_A+\macheps_h\Bigr)
\lVert A\rVert_2.
\qedhere
\end{align*}
\end{proof}

We model the inexact application of the preconditioner \(T\) to $\hat{r}_i$ in
the iteration (\ref{eq:iterate-pinvit-2}) with the equation
\begin{equation}
\label{eq:wm}
\hat{w}_i = T_E\hat{r}_i,
\end{equation}
where \(T_E\) depends on the choice of preconditioner \(T\) and the way to
compute \(T\hat{r}_i\).
Note that \(T_E\) also depends on \(i\).


\begin{theorem}
\label{thm:round-er-pinvit}
Consider the setting of Lemma~\ref{lem:r} and~\eqref{eq:wm}.
If \(\lambda_1<\rho(x_i)<\lambda_2\) and
\[
\gamma := \lVert I-A^{1/2}T_EA^{1/2}\rVert_2
+\gamma_2^h \lVert T_E\rVert_2 \lVert A\rVert_2
+\beta(x_i) \bigl(\macheps_h + (1+\gamma_2^h)\epsilon_r\lVert T_E\rVert_2
\lVert A\rVert_2\bigr)
<1,
\]
with
\[
\beta(x_i)
=\max\biggl\lbrace\frac{\sqrt{\lambda_1\lambda_n}}{\rho(x_i)-\lambda_1},
\frac{\sqrt{\lambda_2\lambda_n}}{\lambda_2-\rho(x_i)}\biggr\rbrace,
\]
then the computed result \(\hat{x}_{i+1}\) of the PINVIT
iteration~\eqref{eq:iterate-pinvit-2} satisfies
\[
\frac{\rho(\hat{x}_{i+1})-\lambda_1}{\lambda_2-\rho(\hat{x}_{i+1})}
\leq \Bigl(\gamma+(1-\gamma)\frac{\lambda_1}{\lambda_2}\Bigr)^2
\frac{\rho(x_i)-\lambda_1}{\lambda_2-\rho(x_i)}.
\]
\end{theorem}
\begin{proof}
By~\eqref{eq:iterate-pinvit-2}, \eqref{eq:wm}, and Lemma~\ref{lem:r}, there
exists a diagonal matrix \(E_0\) (coming from the vector addition) such that
\(\lVert E_0\rVert\leq\macheps_h\) and
\begin{align*}
\label{eq:xm1}
\hat{x}_{i+1}
&=(I+E_0)\bigl(x_i - T_E(I+E)(Ax_i - \rho(x_i)x_i + F x_i)\bigr) \\
&=x_i-\bigl(\tilde{T}_E (A-\rho(x_i) I)-E_0+\tilde{T}_EF\bigr)x_i,
\end{align*}
where \(\tilde{T}_E = (I+E_0)T_E(I+E)\).
Setting $A_\rho = A-\rho(x_i) I$ and using that $x_i = A_\rho^{-1}r_i$,
it follows that
\[
\hat{x}_{i+1}
=x_i-\bigl(\tilde{T}_E-E_0A_\rho^{-1}+\tilde{T}_EFA_\rho^{-1}\bigr)r_i,
\]
which takes the form of PINVIT with a perturbed preconditioner.
This allows us to apply~\cite[Theorem~2.1]{AKNOZ2017}, which requires the
preconditioner to satisfy
\begin{equation}
\label{eq:condpert}
\bigl\lVert I - A^{1/2}(\tilde{T}_E-E_0A_\rho^{-1}+\tilde{T}_EFA_\rho^{-1})
A^{1/2}\bigr\rVert_2
<1.
\end{equation}
We now treat the different terms involved in~\eqref{eq:condpert} separately.
First, we have
\begin{align}
\lVert I-A^{1/2}\Tilde{T}_EA^{1/2}\rVert_2
&\leq \lVert I - A^{1/2} T_E A^{1/2} \rVert_2
+ \lVert A^{1/2}(\Tilde{T}_E - T_E)A^{1/2}\rVert_2 \nonumber \\
&\leq \lVert I - A^{1/2} T_E A^{1/2} \rVert_2
+ \gamma_2^h \lVert T_E\rVert_2 \lVert A\rVert_2.
\label{eq:lala}
\end{align}
By the assumptions, the spectral radius of \(A_\rho^{-1} A^{1/2}\) is given by
\[
\max\biggl\lbrace \frac{\sqrt{\lambda_1}}{\rho(x_i)-\lambda_1},
\frac{\sqrt{\lambda_2}}{\lambda_2-\rho(x_i)} \biggr\rbrace.
\]
This allows us to bound the other terms in~\eqref{eq:condpert} as follows:
\ifnarrowmode
\begin{align*}
&\lVert A^{1/2}(-E_0A_\rho^{-1} + \tilde{T}_EF A_\rho^{-1})A^{1/2}\rVert_2\\
\leq{}&\lVert A^{1/2}\rVert_2 \lVert A_\rho^{-1} A^{1/2}\rVert_2
\bigl(\lVert E_0\rVert_2 + \lVert\tilde{T}_E\rVert_2\lVert F\rVert_2\bigr) \\
\leq{}&\lVert A^{1/2}\rVert_2 \lVert A_\rho^{-1} A^{1/2}\rVert_2
\bigl(\macheps_h + (1+\gamma_2^h)\epsilon_r\lVert T_E\rVert_2\lVert A\rVert_2\bigr) \\
={}&\beta(x_i)\bigl(\macheps_h
+(1+\gamma_2^h)\epsilon_r\lVert T_E\rVert_2\lVert A\rVert_2\bigr).
\end{align*}
\else
\begin{align*}
\lVert A^{1/2}(-E_0A_\rho^{-1} + \tilde{T}_EF A_\rho^{-1})A^{1/2}\rVert_2
&\leq\lVert A^{1/2}\rVert_2 \lVert A_\rho^{-1} A^{1/2}\rVert_2
\bigl(\lVert E_0\rVert_2 + \lVert\tilde{T}_E\rVert_2\lVert F\rVert_2\bigr) \\
&\leq\lVert A^{1/2}\rVert_2 \lVert A_\rho^{-1} A^{1/2}\rVert_2
\bigl(\macheps_h + (1+\gamma_2^h)\epsilon_r\lVert T_E\rVert_2\lVert A\rVert_2\bigr) \\
&=\beta(x_i)\bigl(\macheps_h
+(1+\gamma_2^h)\epsilon_r\lVert T_E\rVert_2\lVert A\rVert_2\bigr).
\end{align*}
\fi
Together with~\eqref{eq:lala}, this implies that the left-hand side
of~\eqref{eq:condpert} is bounded by $\gamma<1$ and the statement of the
theorem follows from~\cite[Theorem~2.1]{AKNOZ2017}.
\end{proof}

\begin{remark}
We remark that the conclusion of Theorem~\ref{thm:round-er-pinvit} does
\emph{not} imply that \(\rho(x_i)-\lambda_1\) can eventually drop below 
machine precision.
For the relative error
\(\bigl(\rho({x}_{i})-\lambda_1\bigr)/\bigl(\lambda_2-\rho({x}_{i})\bigr)\)
to be reduced by the factor
\[
\Bigl(\gamma+(1-\gamma)\frac{\lambda_1}{\lambda_2}\Bigr)^2
=\biggl(\frac{\lambda_1}{\lambda_2}
+\Bigl(1-\frac{\lambda_1}{\lambda_2}\Bigr)\gamma\biggr)^2
\]
during the $i$th iteration, Theorem~\ref{thm:round-er-pinvit} requires that
\[
\lambda_1
<\rho(x_i)
<\lambda_2-\frac{\macheps_h
+(1+\gamma_2^h)\epsilon_r\lVert T_E\rVert_2\lVert A\rVert_2}
{1-\big\lVert I-A^{1/2}T_EA^{1/2}\big\rVert_2
-\gamma_2^h\lVert T_E\rVert_2\lVert A\rVert_2}
\sqrt{\lambda_2\lambda_n}
\]
holds.
This reduction takes place until a Rayleigh quotient \(\rho(\hat{x})\) for an
iterate \(\hat{x}\) is produced for which
\[
\frac{\rho(\hat{x})-\lambda_1}{\sqrt{\lambda_1\lambda_n}}
\leq\frac{\macheps_h+(1+\gamma_2^h)\epsilon_r\lVert T_E\rVert_2\lVert A\rVert_2}
{1-\lVert I-A^{1/2}T_EA^{1/2}\rVert_2
-\gamma_2^h\lVert T_E\rVert_2\lVert A\rVert_2}.
\]
For reasonable choices of $T_E$, this means that the error is reduced until it
reaches the level of working precision.
\end{remark}


The quantity \(\lVert I - A^{1/2}T_EA^{1/2}\rVert_2\) critically determines the convergence rate of PINVIT. The following lemma provides an estimate if $T_E$ corresponds to applying $A^{-1}$ in low precision via the Cholesky factorization.
\begin{lemma}
\label{lem:T}
Suppose that the application of the preconditioner $T$ in one step of
PINVIT~\eqref{eq:iterate-pinvit-2} is implemented by applying $A^{-1}$ in low
precision, via performing the Cholesky factorization of $A$ followed by
forward and backward substitution.
If \(\epsilon_T:=4n(3n+1)\kappa(A)\macheps_l<1\), where
\(\kappa(A) = \lVert A\rVert_2\lVert A^{-1}\rVert_2\) and $\macheps_l$ denotes
unit roundoff in low precision, then
\[
\bigl\lVert I-A^{1/2}T_EA^{1/2}\bigr\rVert_2
\leq\frac{\epsilon_T}{1-\epsilon_T}.
\]
\end{lemma}

\begin{proof}
Using~\cite[Theorem~10.4]{H2002}, there exists a symmetric matrix \(E_0\) such
that
\[
\hat{w}_i=\fl(T\hat{r}_i)=(A+E_0)^{-1}\hat{r}_i, \qquad
\lVert E_0\rVert_2\leq 4n(3n+1)\macheps_l\lVert A\rVert_2,
\]
which means \(T_E=(A+E_0)^{-1}\) and, moreover,
\[
A^{1/2}T_EA^{1/2}=\bigl(I+A^{-1/2}E_0A^{-1/2}\bigr)^{-1}.
\]
Then by
\(\bigl\lVert A^{-1/2}E_0A^{-1/2}\bigr\rVert_2%
\leq 4n(3n+1)\kappa(A)\macheps_l<1\),
we have
\[
\bigl(I+A^{-1/2}E_0A^{-1/2}\bigr)^{-1}
=\sum_{i=0}^{\infty}\bigl(-A^{-1/2}E_0A^{-1/2}\bigr)^i.
\]
Thus, it holds that
\ifnarrowmode
\begin{align*}
\bigl\lVert I-A^{1/2}T_EA^{1/2}\bigr\rVert_2
&=\bigl\lVert I-\bigl(I+A^{-1/2}E_0A^{-1/2}\bigr)^{-1}\bigr\rVert_2\\
&\leq\sum_{i=1}^{\infty}\bigl\lVert A^{-1/2}E_0A^{-1/2}\bigr\rVert_2^i\\
&\leq\frac{\epsilon_T}{1-\epsilon_T}.
\qedhere
\end{align*}
\else
\[
\bigl\lVert I-A^{1/2}T_EA^{1/2}\bigr\rVert_2
=\bigl\lVert I-\bigl(I+A^{-1/2}E_0A^{-1/2}\bigr)^{-1}\bigr\rVert_2
\leq\sum_{i=1}^{\infty}\bigl\lVert A^{-1/2}E_0A^{-1/2}\bigr\rVert_2^i
\leq\frac{\epsilon_T}{1-\epsilon_T}.
\qedhere
\]
\fi
\end{proof}

\section{Numerical experiments}
\label{sec:experiments}
In this section, we present numerical results for our mixed
precision LOBPCG algorithm.
In our tests, the working precision is IEEE double precision and the lower
precision is IEEE single precision.
Most tests are performed on a Linux server equipped with two twelve-core Intel
Xeon E5-2670~v3 2.30~GHz CPUs and two Nvidia GeForce GTX~1080 GPUs.
The tests in Section~\ref{subsec:A30} also use an Nvidia A30 GPU.
There are 128~GB of main memory on the CPUs and 11,178.6~MB of main memory on
each GPU.
Our program uses only one GPU and one thread on the CPU.

\subsection{Experiment settings}
In our experiments we compute a few smallest eigenvalues and the corresponding
eigenvectors of Hermitian matrices using the LOBPCG algorithm.
The following variants of the LOBPCG algorithm are tested:
\begin{enumerate}
\item
{\bf DLOBPCG-dchol}: LOBPCG algorithm performed entirely in double precision.
\item
{\bf DLOBPCG-schol}: LOBPCG algorithm performed in double precision, except for
single precision preconditioning.
\item
{\bf MPLOBPCG-schol}: mixed precision LOBPCG algorithm
(Algorithm~\ref{alg:mlobpcg}) with single precision preconditioning and
initial guess computed by the single precision LOBPCG algorithm;
mixed precision orthogonalization (Algorithm~\ref{alg:mqr}) is also used.
\end{enumerate}

When to computing \(k\) eigenpairs, we run LOBPCG algorithm with a block size that is about
\(50\%\) larger in order to enhance robustness.
The algorithm terminates once the \(k\) smallest eigenvalues and the
corresponding eigenvectors converge.
The convergence criterion is
\begin{equation}
\label{eq:stop}
\lVert AX(:, j)-\Theta(j, j)X(:, j)\rVert_2
\leq \tol\cdot(\lVert A\rVert_2+\lvert \Theta(j, j)\rvert)
\lVert X(:, j)\rVert_2,
\end{equation}
where \(\lVert A\rVert_2\) is estimated through \(\lVert A\rVert_2\approx%
\lVert\Omega A\rVert_{\fro}/\lVert\Omega\rVert_{\fro}\) using a Gaussian
random matrix \(\Omega\in\mathbb C^{m\times n}\) with \(m\ll n\).
The threshold \(\tol\) in~\eqref{eq:stop} is set to \(10^{-12}\) for all these
three algorithms, and is \(5\times10^{-6}\) when computing a good initial
guess for MPLOBPCG-schol.

In our tests, we use \(\Pi\iherm L\iherm L^{-1}\Pi^{-1}\) as the
preconditioner for Algorithm~\ref{alg:mlobpcg}, where \(\Pi\) is a permutation
matrix, and \(L\) is the (pivoted) Cholesky factor of \(A\) satisfying
\(\Pi\herm A\Pi=LL\herm\) computed in single precision.
The preconditioning stage in DLOBPCG-schol/MPLOBPCG-schol is to compute
\(W_\single=\Pi L\iherm L^{-1}\Pi\herm\single(R)\) by solving two triangular
systems in single precision.
In practice, we apply \(\Pi\) to the given matrix \(A\) instead of applying
\(\Pi\) to \(\single(R)\) in each iteration.
We can benefit from it if \(A\) is sufficiently sparse or the convergence of
LOBPCG is not too rapid (i.e., it takes many iterations to converge).

We test the LOBPCG algorithm for both sparse matrices and dense matrices.
Table~\ref{table:implementation} summarizes the software libraries used under
different settings.
The CHOLMOD package~\cite{CDHR2008} can compute sparse Cholesky factorization
on both CPU and GPU, while triangular linear solvers are only supported only
in CPU.
Note that CHOLMOD was developed only for double precision arithmetic; we have
derived a single precision version for the purpose of our tests.

\begin{table}[!tb]
\centering
\caption{Libraries used in our implementation.}
\label{table:implementation}
\begin{tabular}{ccccc}
\hline
& Cholesky & \texttt{TRSM} & mat--vec & others \\
\hline
CPU/sparse & CHOLMOD & CHOLMOD  & MKL      & LAPACK \\
GPU/sparse & CHOLMOD & cuSPARSE & cuSPARSE & MAGMA \\
CPU/dense  & LAPACK  & LAPACK   & LAPACK   & LAPACK \\
GPU/dense  & MAGMA   & MAGMA    & MAGMA    & MAGMA \\
\hline
\end{tabular}
\end{table}

\subsection{Advantage of mixed precision orthogonalization}
Before discussing the LOBPCG algorithm, we first report the run time and
savings of the mixed precision orthgonalization algorithm (i.e.,
Algorithm~\ref{alg:mqr}) in Table~\ref{table:mqr}.
We can see that for tall-skinny matrices Algorithm~\ref{alg:mqr} can reduce
the run time by a factor of \(1/4\)--\(1/3\) compared to \texttt{DGEQRF} in
cuSOLVER.
Thus, it is worth using this mixed approach for orthogonalization.

\begin{table}[!tb]
\centering
\caption{Run time of Algorithm~\ref{alg:mqr} in seconds.}
\label{table:mqr}
\begin{tabular}{cccc}
\hline
\((n,k)\) & Run time of Algorithm~\ref{alg:mqr} &
Run time of \texttt{DGEQRF} & Savings \\
\hline
\((40000,20)\) & \(1.544\times10^{-3}\) & \(2.078\times10^{-3}\) & \(25.7\%\) \\
\((40000,25)\) & \(1.907\times10^{-3}\) & \(2.823\times10^{-3}\) & \(32.5\%\) \\
\((40000,30)\) & \(2.371\times10^{-3}\) & \(3.579\times10^{-3}\) & \(33.7\%\) \\
\((40000,35)\) & \(3.586\times10^{-3}\) & \(4.886\times10^{-3}\) & \(26.6\%\) \\
\((40000,40)\) & \(4.000\times10^{-3}\) & \(5.726\times10^{-3}\) & \(30.1\%\) \\
\((40000,45)\) & \(4.495\times10^{-3}\) & \(6.862\times10^{-3}\) & \(34.5\%\) \\
\hline
\end{tabular}
\end{table}

\subsection{Tests for sparse matrices}
We choose six sparse positive definite matrices from from the SuiteSparse
Matrix Collection.%
\footnote{\url{https://sparse.tamu.edu}}
Table~\ref{table:info_of_matrices} shows the information of these
sparse matrices.
We compute \(30\) eigenpairs using a randomly generated initial guess with
\(45\) columns for each
matrix, and report the relative run time, which is the ratio of the wall clock
time of a solver over the wall clock time of DLOBPCG-dchol.

\begin{table}[!tb]
\centering
\caption{Information of sparse testing matrices.}
\label{table:info_of_matrices}
\begin{tabular}{ccccc}
\hline
Name& Size& NNZ& Sparsity& NNZ of \(L\)\\
\hline
obstclae        & \hphantom{0}40,000 & \hphantom{0,}197,608 & \(1.235\times10^{-4}\)& \hphantom{}1,561,880\\
shallow\_water2 & \hphantom{0}81,920 & \hphantom{0,}327,680 & \(4.883\times10^{-5}\)& \hphantom{}3,483,014\\
Dubcova2        & \hphantom{0}65,025 & \hphantom{}1,030,225 & \(2.437\times10^{-4}\)& \hphantom{}3,804,558\\
Dubcova3        & \hphantom{}146,689 & \hphantom{}3,636,643 & \(1.690\times10^{-4}\)& \hphantom{}7,409,077\\
finan512        & \hphantom{0}74,752 & \hphantom{0,}596,992 & \(1.068\times10^{-4}\)& \hphantom{}3,376,835\\
2D-Laplace      & \hphantom{0}25,000 & \hphantom{0,}114,990 & \(1.840\times10^{-4}\)& \hphantom{0,}466,491\\
\hline
\end{tabular}
\end{table}

\begin{figure}[!tb]
\centering
\includegraphics[width=1.0\textwidth]{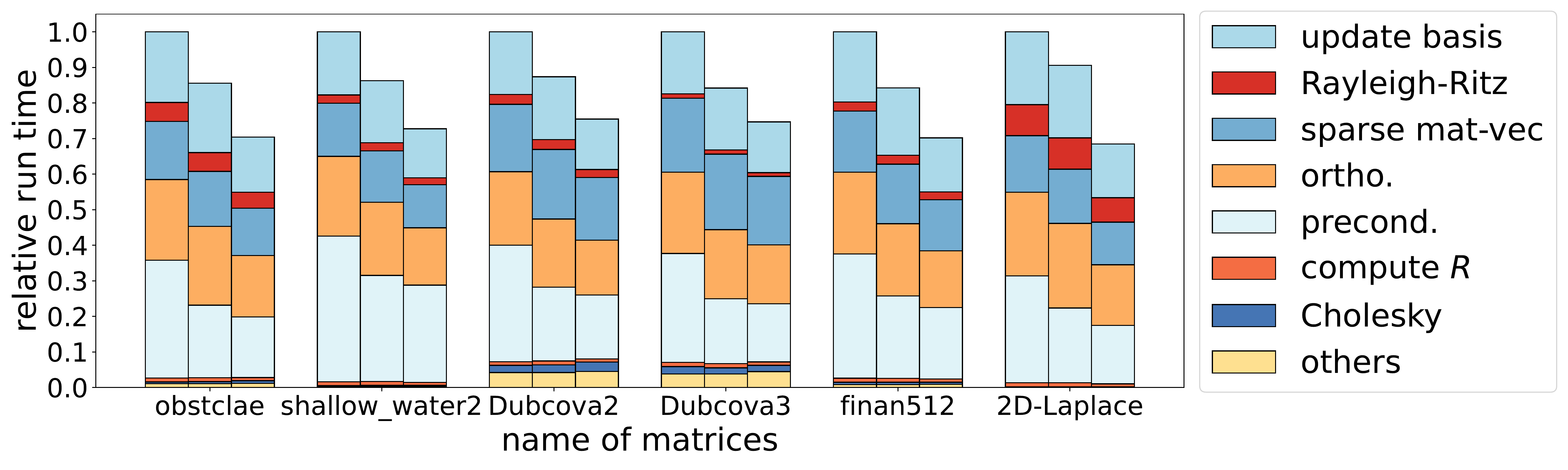}
\caption{Tests for real sparse matrices on CPU.
For each matrix, the three columns from left to right represent the result of
DLOBPCG-dchol, DLOBPCG-schol, and MPLOBPCG-schol, respectively.}
\label{fig:cpu_sparse}
\end{figure}

\begin{figure}[!tb]
\centering
\includegraphics[width=1.0\textwidth]{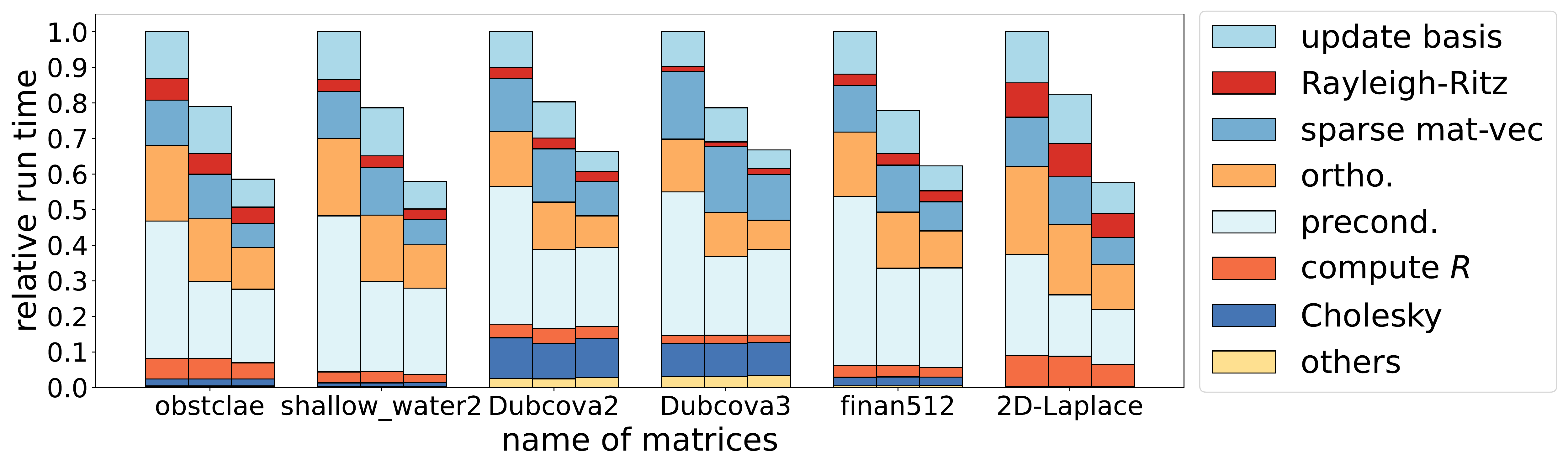}
\caption{Tests for real sparse matrices on GPU.
For each matrix, the three columns from left to right represent the result of
DLOBPCG-dchol, DLOBPCG-schol, and MPLOBPCG-schol, respectively.}
\label{fig:gpu_sparse}
\end{figure}

Figures~\ref{fig:cpu_sparse} and~\ref{fig:gpu_sparse}, respectively, show the
relative run time on CPU and GPU.
For all test cases, preconditioning in single precision reduces the execution
time of the LOBPCG algorithm.
Using an initial guess computed by the single precision LOBPCG algorithm and
adopting mixed precision orthogonalization makes the algorithm more efficient.
Compared to DLOBPCG-dchol, MPLOBPCG-schol is about \(1.43\times\) faster on
CPU, and is about \(1.67\times\) faster on GPU.

We should also mention that the number of iterations for different variants
of the LOBPCG algorithms are similar, though they are not shown in the
figures.
Sometimes MPLOBPCG-schol can require fewer iterations to converge because
there is a restart when we use the lower precision result as the initial
guess.
For instance, the total iterations of DLOBPCG-dchol, DLOBPCG-schol and
MPLOBPCG-schol are \(533\), \(534\), and \(461\), respectively, for the
2D-Laplace matrix.

\subsection{Tests for dense matrices}
We also test the LOBPCG algorithm for a few dense matrices which are popular
in machine learning.
These dense matrices are kernel matrices generated by certain kernel
functions as follows.
Let \(x_1\), \(x_2\), \(\dotsc\), \(x_n\in\mathbb R^n\) be uniform random
vectors generated by \(\mathtt{XLARNV}\) from LAPACK.
We construct a matrix \(K\) by applying the Gaussian kernel function
\[
K_{ij}=k(x_i,x_j)=\me^{-\lVert x_i-x_j\rVert_2/2}.
\]
Similarly, we can apply the polynomial kernel function
\[
k(x_i,x_j)=(x_i\trans x_j+1)^3
\]
to construct another kernel matrix.
Using two sets of random vectors \(\lbrace{x_1,x_2,\dotsc,x_n\rbrace}\) and
\(\lbrace{y_1,y_2,\dotsc,y_n\rbrace}\) in \(\mathbb R^n\), we also construct
complex kernel matrices through
\[
K_{ij}=k(x_i,x_j)+k(y_i,y_j)+\mi\bigl(k(x_i,y_j)-k(y_i,x_j)\bigr),
\]
where \(k(\cdot,\cdot)\) is either the Gaussian kernel function
or the polynomial kernel function.

\begin{figure}[!tb]
\centering
\includegraphics[width=1.0\textwidth]{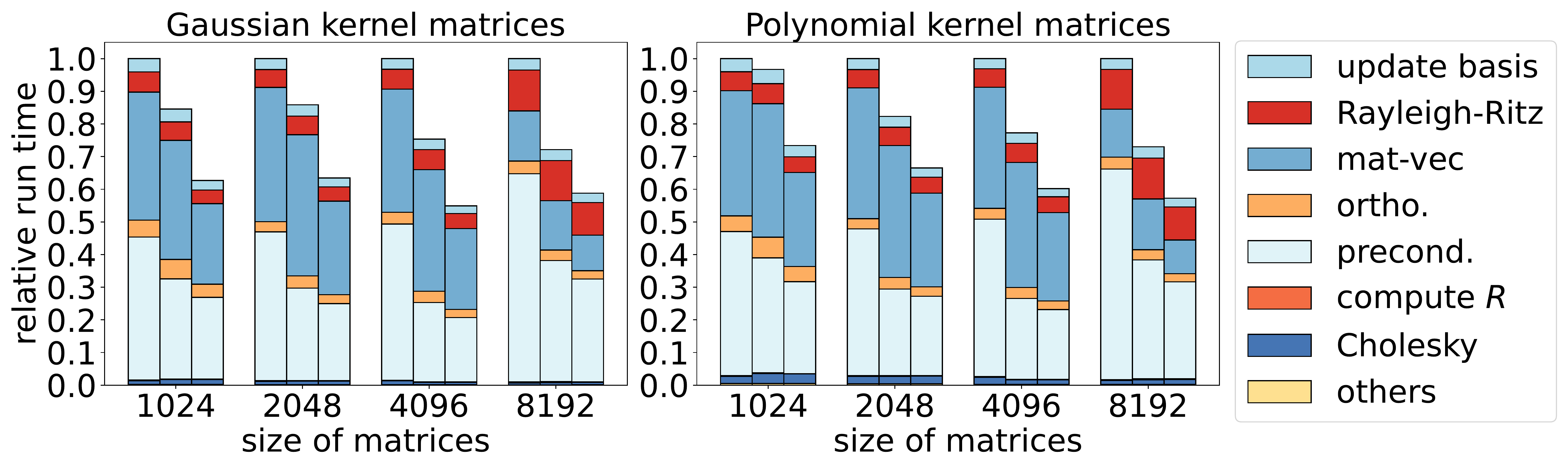}
\caption{Tests for real dense kernel matrices on CPU.
For each matrix, the three columns from left to right represent the result of
DLOBPCG-dchol, DLOBPCG-schol, and MPLOBPCG-schol, respectively.}
\label{fig:cpu_dense}
\end{figure}

\begin{figure}[!tb]
\centering
\includegraphics[width=1.0\textwidth]{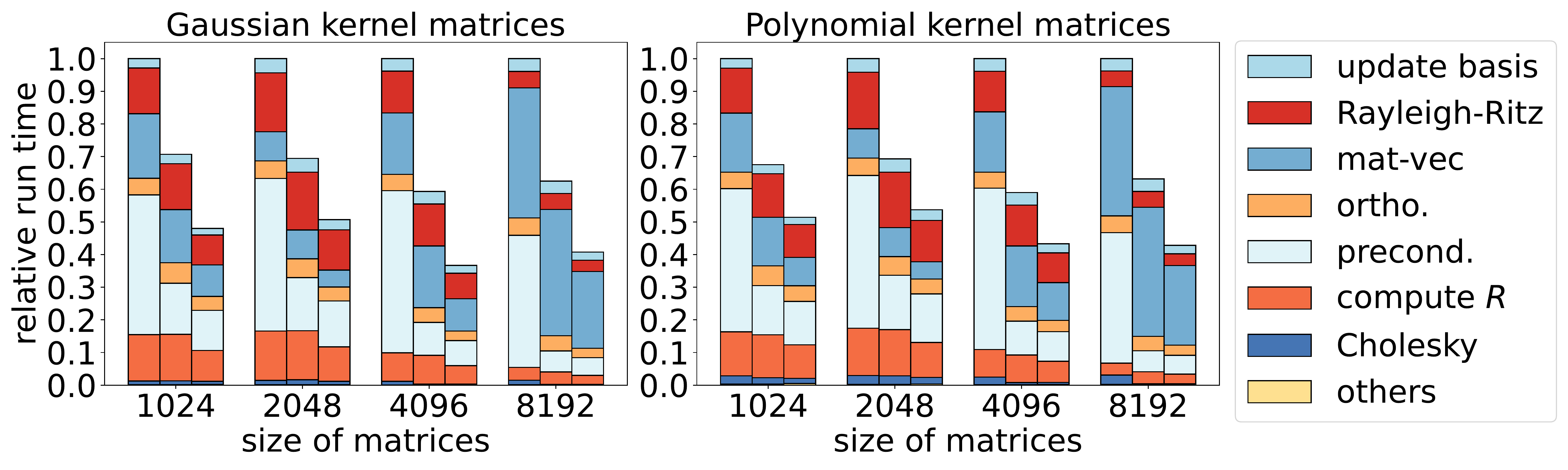}
\caption{Tests for real dense kernel matrices on GPU.
For each matrix, the three columns from left to right represent the result of
DLOBPCG-dchol, DLOBPCG-schol, and MPLOBPCG-schol, respectively.}
\label{fig:gpu_dense}
\end{figure}

\begin{figure}[!tb]
\centering
\includegraphics[width=1.0\textwidth]{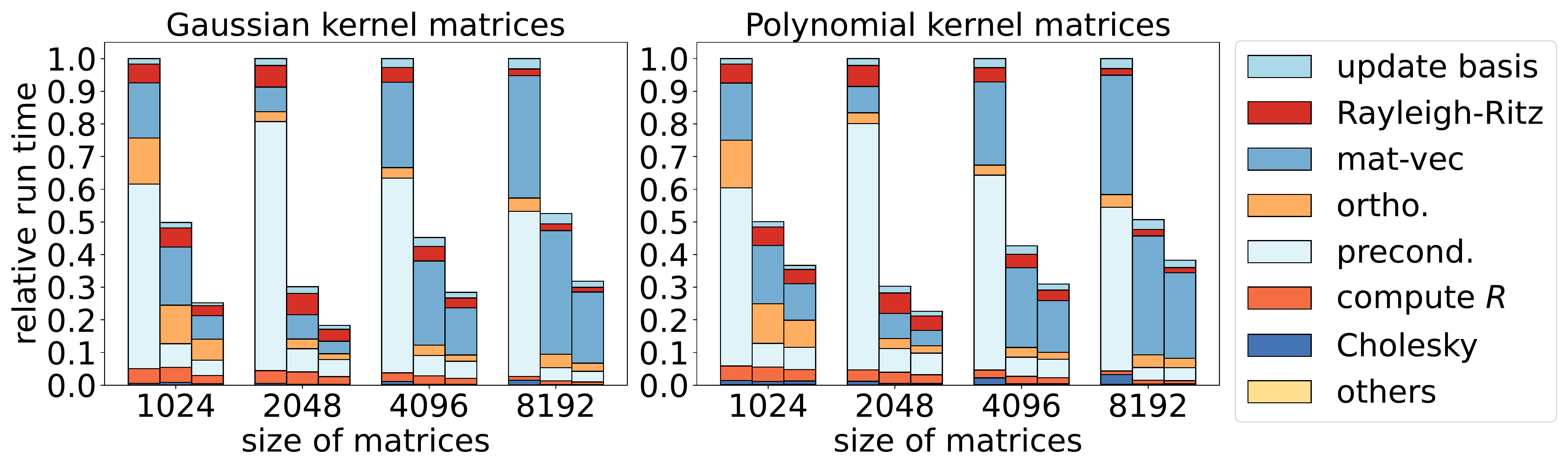}
\caption{Tests for complex dense kernel matrices on GPU.
For each matrix, the three columns from left to right represent the result of
DLOBPCG-dchol, DLOBPCG-schol, and MPLOBPCG-schol, respectively.}
\label{fig:gpu_dense_complex}
\end{figure}

We choose \(n\in\lbrace{1024,2048,4096,8192\rbrace}\) in our experiments, and
compute \(5n/1024\) smallest eigenvalues and the corresponding eigenvectors.
The rank of initial guess is chosen as \(8n/1024\) accordingly.
Figures~\ref{fig:cpu_dense}, \ref{fig:gpu_dense},
and~\ref{fig:gpu_dense_complex} show the relative run time of different
variants of the LOBPCG algorithm.
For real matrices, MPLOBPCG-schol achieves \(1.67\times\) and \(2\times\)
speedup compared to DLOBPCG-dchol on CPU and GPU, respectively.
The speedup is higher than that for sparse matrices, because dense matrices
are more compute-intensive.
The benefit for mixed precision approaches is more significant for complex
matrices---the speedup becomes over \(2.5\times\) and up to
\(5\times\) on GPU.

\subsection{Tests on different GPUs}
\label{subsec:A30}
By far our tests are performed with an Nvidia GeForce GTX~1080 GPU, which is a
consumer-grade GPU.
In fact, there are two different types of GPU---consumer-grade and
server-grade.
Compared to consumer-grader GPUs server-grade GPUs usually have better
hardware support for double precision arithmetic.
Hence the performance difference between single and double precision
arithmetic is larger on consumer-grade GPUs.

\begin{figure}[!tb]
\centering
\includegraphics[width=0.8\textwidth]{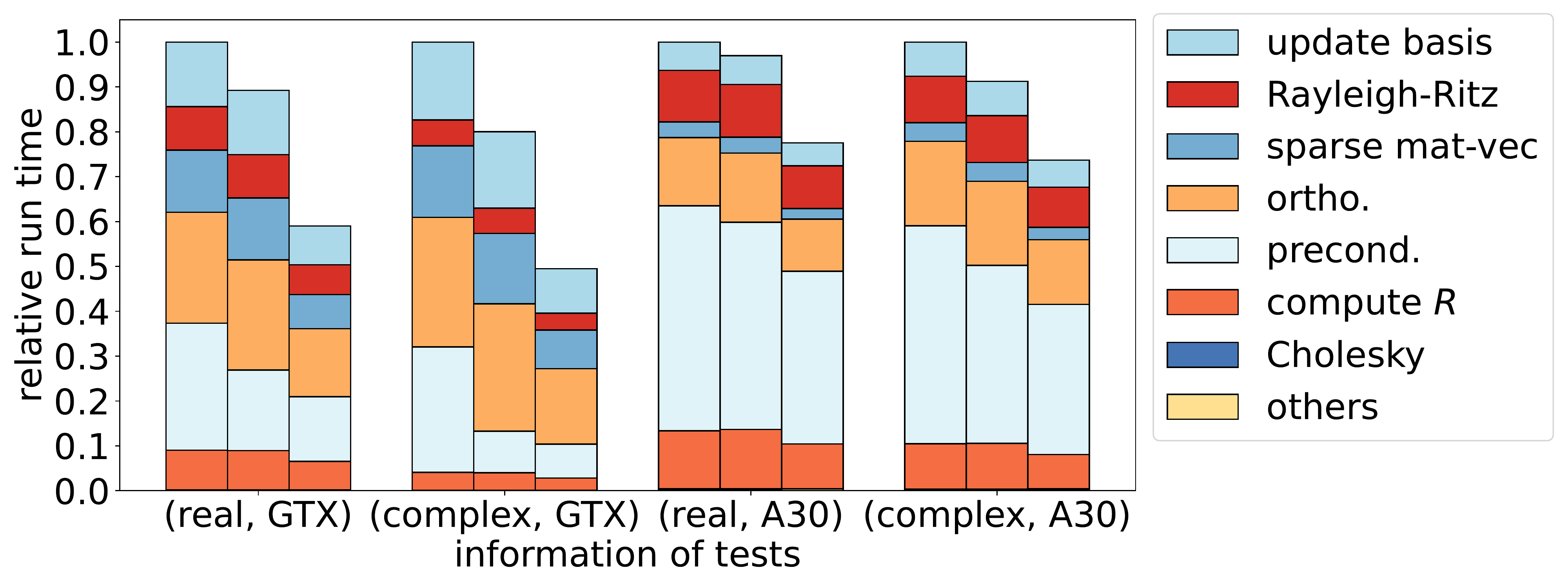}
\caption{Tests for real and complex 2D-Laplace matrices in different GPUs.
For each case, the three columns from left to right represent the result of
DLOBPCG-dchol, DLOBPCG-schol, and MPLOBPCG-schol, respectively.}
\label{fig:laplace}
\end{figure}

In the following we report some results collected from runs on an Nvidia A30
GPU, which is a server-grade one.
We use the matrix 2D-Laplace in this test.
By perturbing off-diagonal entries of this matrix by \(\pm10^{-16}\cdot\mi\),
we also obtain a Hermtian positive definite matrix for testing complex
arithmetic.
From Figure~\ref{fig:laplace}, it can be seen that single precision has
limited advantage over double precision on this server-grade GPU.
Though MPLOBPCG-schol still achieves about \(1.3\times\) speedup compared to
DLOBPCG-dchol, the benefit for adopting single precision arithmetic is much
lower than that on NVIDIA GeForce GTX-1080 which is a consumer-grade GPU.

\section{Conclusion}
\label{sec:conclusion}
In this paper, we have proposed a mixed precision LOBPCG algorithm with a preconditioner based on a (sparse)
Cholesky factorization.
Both the initial guess and the preconditioner are computed in reduced
precision.
This largely improves the performance while it only has marginal impact on 
convergence.
In our mixed precision LOBPCG algorithm, orthogonalization is also performed
in a mixed precision manner to further improve performance.
We analyze the rounding error of the PINVIT algorithm, which can be viewed as
a simplified version of the LOBPCG algorithm, to confirm that our mixed
precision algorithm is as accurate as the fixed precision one.
Numerical experiments illustrate that adopting mixed precision arithmetic can
significantly accelerate the execution of the LOBPCG algorithm on both CPUs
and GPUs.

\backmatter

\bmhead{Acknowledgments}
The authors thank Erin Carson for helpful discussions.
Part of this work was performed when the second author was visiting EPF
Lausanne in 2022.

Yuxin Ma is partially supported by the State Scholarship Fund of China Scholarship Council (CSC) under Grant No.~202106100093, National Key R\&D Program of China under
Grant No.~2021YFA1003305 and National Natural Science Foundation of China
under Grant No.~71991471.
Meiyue Shao is partially supported by by the National Natural Science
Foundation of China under grant No.~11971118.

\bibliography{mybib}

\end{document}